\documentclass[review]{elsarticle}

\usepackage{amssymb,   amsmath,   amssymb,   lineno,  subfigure}
\modulolinenumbers[5]

\usepackage{graphicx}

\journal{Journal of Combinatorial Theory,   Series A}

\begin{document}

\begin{frontmatter}

\title{A New Result on Packing Unit Squares into a Large Square}

\author[mymainaddress]{Shuang Wang}
\ead{wangshuang@jlu.edu.cn}

\author[mymainaddress]{Tian Dong\corref{mycorrespondingauthor}\fnref{myfootnote}}
\cortext[mycorrespondingauthor]{Corresponding author}
\ead{dongtian@jlu.edu.cn}
\fntext[myfootnote]{Tian Dong was supported by National Natural Science Foundation of China under Grant No. 11101185 and 11171133.}

\author[mymainaddress]{Jiamin Li}
\ead{jmli@jlu.edu.cn}

\address[mymainaddress]{School of Mathematics,   Jilin University,   Changchun,   Jilin 130012,   China}

\begin{abstract}
In their 2009 note: \emph{Packing equal squares into a large square},  Chung and Graham proved that the wasted area of a large square of side length $x$ is $O\left(x^{(3+\sqrt{2})/7}\log x\right)$ after maximum number of non-overlapping unit squares are packed into it,  which improved the earlier results of Erd\H{o}s-Graham and Karabash-Soifer. Here we further improve the result to $O(x^{5/8})$ that also leads to an improvement of the bound for the dual problem: finding the minimum number of unit squares needed for covering the large square,  from $x^2+O\left(x^{(3+\sqrt{2})/7}\log x\right)$ to $x^2+O(x^{5/8})$.
\end{abstract}


\newtheorem{thm}{Theorem}
\newproof{pf}{Proof}
\newdefinition{defn}{Definition}
\newdefinition{rmk}{Remark}

\end{frontmatter}

\linenumbers

\section{Introduction}
\label{sec:int}

In 1975,  Erd\H{o}s and Graham \cite{EG1975} investigated the problem of packing a square of side length $x$ with as many non-overlapping unit squares as possible. In other words, the wasted area should be as small as possible. From then on, the problem have already been well studied in the literature \cite{RV1978, Str1984, Str2003, KS2008, CG2009, Fri2009, Ben2010}, in which \cite{RV1978, KS2008, CG2009} focus on the case when $x$ is large enough. Following \cite{KS2008}, we call the problem Packing Waste Problem. Also, there is a dual problem, called Covering Waste Problem in \cite{KS2008}, which is concerned with covering the square with minimium number of unit squares\cite{KS2008,CG2009,KS2006, Soi2006, FP2009, Jan2009}.

Erd\H{o}s and Graham obtained the first estimation of Packing Waste Problem as $O(x^{7/11})$ \cite{EG1975}. Later, D. Karabash and A. Soifer in \cite{KS2006} gave the estimation of Covering Waste Problem as $O(x^{2/3})$ that was improved in \cite{KS2008} to $O(x^{7/11})$. In 2009, Chung and Graham \cite{CG2009} found the best previous bound $O\left(x^{(3+\sqrt{2})/7}\log x\right)$ for both problems.

In this paper we use basic analysis tools to improve the result of Chung and Graham to $O(x^{5/8})$ also for both problems.

\section{Preliminary}\label{pre}


Let $A$ be a closed planar region and $S(A)$ the area of it. We define two functions
\begin{align*}
  W(A)&=S(A)-\sup s(A_\lambda),\\
  W'(A)&=\inf s(A'_\lambda)-S(A),
\end{align*}
where $A_\lambda \subset A$ is a union set of non-overlapping unit squares, and $A'_\lambda \supset A$ is a union set of unit squares (non-overlapping is not necessary). Specially, when $A$ is a square of side length $x$,  we denote $W(A), W'(A)$ as $W(x), W'(x)$ respectively.

\begin{figure}[!htph]
\centering
\includegraphics[width=0.4\textwidth]{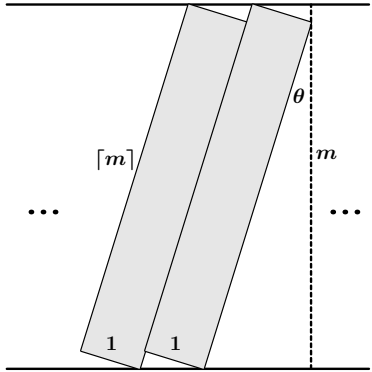}
\caption{Packing a strip of width $m$.}\label{f:1}
\end{figure}

To our opinion, the basic task of Packing or Covering Waste Problem is packing or covering a strip of non-integer width \cite{CG2009}, say $m$. Basic idea for packing a strip \cite{CG2009} is to pack stacks of non-overlapping unit squares of height $\lceil m \rceil$ into the strip as close to being orthogonal as possible (see Fig. \ref{f:1}), namely minimize the angle $\theta$ in Fig. \ref{f:1} which satisfies
\begin{equation}\label{e:1}
  \lceil m \rceil \cos \theta+\sin\theta=m.
\end{equation}
Let $r=m-\lfloor m\rfloor$. Obviously when $r=0$, $\theta=0$ trivially. Otherwise, we let $\theta=\alpha m^\beta+o(m^\beta)$. By comparing with the constant term of \eqref{e:1}, we have
$$
\theta=\sqrt{2-r}~m^{-1/2}+o(m^{-1/2}).
$$

\begin{figure}[!htph]
\centering
\includegraphics[width=0.4\textwidth]{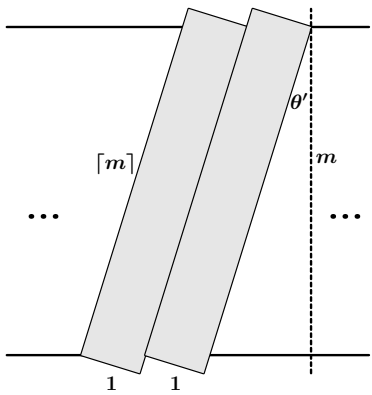}
\caption{Covering a strip of width $m$.}\label{f:2}
\end{figure}

Similarly, as shown in Fig. \ref{f:2}, we also use stacks of unit squares of height $\lceil m \rceil$ (hereafter we will call the stacks as rectangles of size $ 1 \times \lceil m \rceil $ for simplicity) to cover the strip, then angle $\theta'$ in Fig. \ref{f:2} satisfies
\begin{equation}\label{e:2}
\lceil m \rceil \cos \theta'-\sin\theta'=m.
\end{equation}
We also have $\theta'=0$ when $r=0$. If not, then
$$
\theta'=\sqrt{2-r}~m^{-1/2}+o(m^{-1/2}).
$$
Note that when $m\rightarrow \infty, \theta$ and $\theta'$ are less than $\sqrt{2}~m^{-1/2}$.

\section{Packing Waste Problem}\label{PWP}

In this section, we will present our main result on Packing Waste Problem in Theorem \ref{t:1}. For the proof of it,  three types of basic shapes are introduced as follows.

\begin{figure}[!htph]
\centering
\subfigure[Type 1 shape.]{\includegraphics[width=0.3\textwidth ]{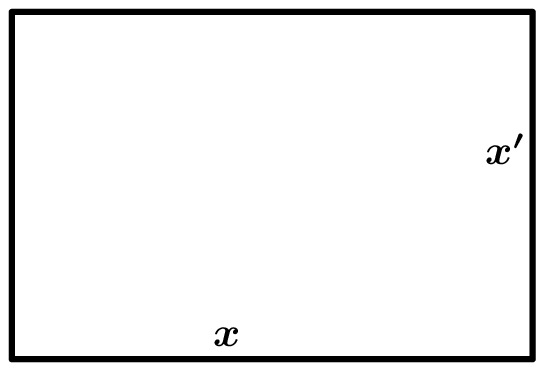}\label{sf:1}}
\subfigure[Type 2 shape.]{\includegraphics[width=0.3\textwidth]{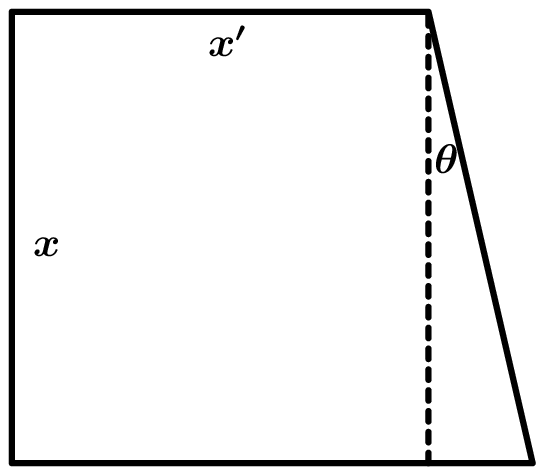}\label{sf:2}}
\subfigure[Type 3 shape.]{\includegraphics[width=0.3\textwidth]{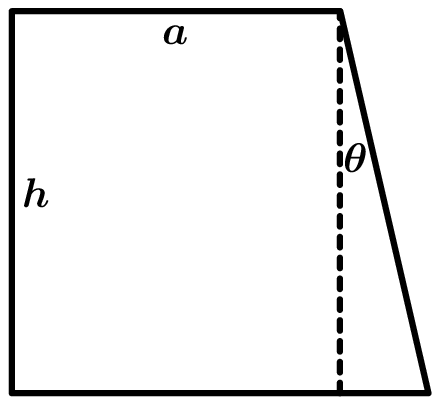}\label{sf:3}}
\caption{Three types of basic shapes.}\label{f:3}
\vskip 0.3cm
\end{figure}

\begin{description}
  \item[Type 1 shape] Rectangle $T_1$  has a length $x$ and width $x'$ (see subfigure \subref{sf:1} of Fig. \ref{f:3}) satisfying $ x^{3/4}\leq x'\leq c x$ with $c\leq 7$ a constant.
  \item[Type 2 shape] Trapezoid $T_2$ has a height of $x$,   a top edge of length $x'$  (see subfigure \subref{sf:2} of Fig. \ref{f:3}) satisfying $x'\sim 2x^{1/2}$ and the angle $\theta$ between the right-hand side and a vertical line satisfying $0 < \theta<\sqrt{2}x^{-1/2}.$
  \item[Type 3 shape] Trapezoid $T_3$ has a height $h \sim \frac{1}{2}x^{1/2}$ and a top edge of length $a$ (see subfigure \subref{sf:3} of Fig. \ref{f:3}) where $a=\lfloor x^{1/3}+\sqrt{2}~x^{1/6}\rfloor$ is an exact integer. The angle $\theta$ between the right-hand side and a vertical line satisfies $0<\theta<\sqrt{2}x^{-1/2}$.
\end{description}

The proof of Theorem \ref{t:1} will be completed by an induction based on effective packings of these shapes.

\begin{thm}\label{t:1}
Keep the notations above. Then
\begin{enumerate}[(i)]
  \item $W(T_1)\leq ((15+c)\sqrt{2}+38)x^{5/8}$.
  \item $W(T_2)\leq (\frac{19}{2}+\frac{7}{2}\sqrt{2})x^{5/6}$.
  \item $W(T_3)\leq (\frac{19}{4}+\frac{7}{4}\sqrt{2})x^{1/3}$.
\end{enumerate}
Specially, when $T_1$ is a square of side length $x$, then $W(x) \leq (16\sqrt{2}+38)x^{5/8}$.
\end{thm}

\begin{pf}
(i) We partition Type 1 rectangle $T_1$ into a rectangle $S_1$ of size $m_1 \times (x-m_2)$, a rectangle $S_2$ of size $m_2 \times x'$, and an integer-sided rectangle $T'_1$, where $m_1, m_2 \sim m=x^{3/4}$, as shown in Fig. \ref{f:4}. It is easy to see that $T'_1$ can be perfectly packed, that is $W(T'_1)=0$. Next, we pack $S_1$ and $S_2$ with rectangles of size $1 \times \lceil m_1 \rceil$ and $1 \times \lceil m_2 \rceil$ respectively. Finally, only four regions $T_{2i}, i=1 , 2, 3, 4$, at each end of $S_1$ and $S_2$, left unfilled which clearly belong to Type 2 with height about $m$, a top edge of length $m'\sim 2m^{1/2}$, and $\theta < \sqrt{2}m^{-1/2}$.

\begin{figure}[!htph]
\centering
\includegraphics[width=0.6\textwidth]{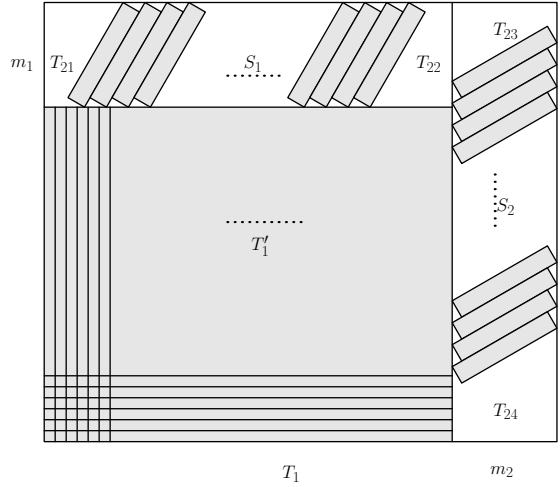}
\caption{Packing Type 1 rectangle.}\label{f:4}
\end{figure}

Applying (ii), the wasted area
\begin{eqnarray*}
W(T_1)&\leq&0 + x\cdot 2 \cdot \frac{1}{2}\tan\theta+x'\cdot 2 \cdot \frac{1}{2}\tan\theta+\sum_{i=1}^4W(T_{2i})\\
      &\leq& (x+x')\cdot \sqrt{2}m^{-1/2}+4(\frac{19}{2}+\frac{7}{2}\sqrt{2})m^{5/6}\\
      &\leq& ((15+c)\sqrt{2}+38)x^{5/8}.
\end{eqnarray*}
Specially, when $T_1$ is a square of side length $x$, $W(x) \leq (16\sqrt{2}+38)x^{5/8}$.

(ii) Now we partition the Type 2 trapezoid $T_2$ into rectangles $A_1, \cdots, A_s$ and Type 3 trapezoids $B_1, \cdots, B_s$ (see Fig. \ref{f:5}). Each $B_i$ has height $h \sim \frac{1}{2}x^{1/2}$ and top edge of length integer $a$. Thus, $s \sim 2x^{1/2}$.

\begin{figure}[!htph]
\centering
\includegraphics[width=0.4\textwidth]{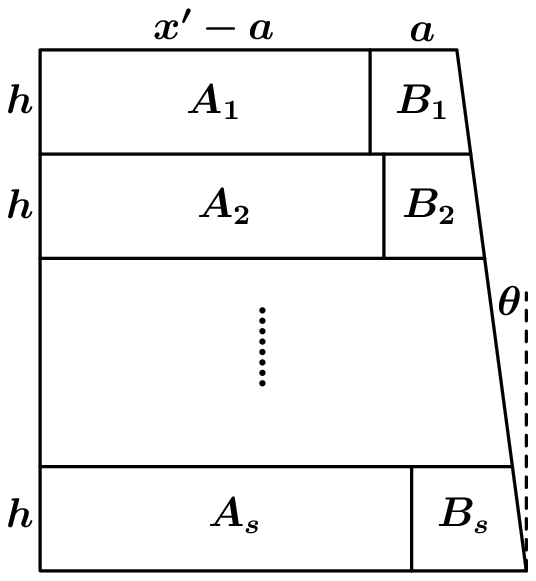}
\caption{Packing Type 2 trapezoid.}\label{f:5}
\end{figure}

Let $a_i$ be the width of $A_i$. Then we have $x^{1/2}<a_i<(2+\sqrt{2})x^{1/2},
2h<a_i<2(2+\sqrt{2})h$. From (i), we obtain $W(A_i)=O(h^{5/8})=O(x^{5/16})$, hence
$$
W\left(\bigcup_{i=1}^{s} A_i\right)\leq \sum_{i=1}^{s} W(A_i)\leq O(x^{5/16})\cdot s=O(x^{13/16}).
$$
Further, (iii) implies that
$$
W\left(\bigcup_{i=1}^{s} B_i\right)\leq \sum_{i=1}^{s} W(B_i)
      \leq \left(\frac{19}{4}+\frac{7}{4}\sqrt{2}\right)x^{1/3}\cdot s
      \leq (\frac{19}{2}+\frac{7}{2}\sqrt{2})x^{5/6},
$$
which leads to the wasted area of $T_2$
$$
W(T_2) \leq W\left(\bigcup_{i=1}^{s} A_i\right)+W\left(\bigcup_{i=1}^{s} B_i\right)\leq \left(\frac{19}{2}+\frac{7}{2}\sqrt{2}\right)x^{5/6}.
$$

(iii) We will partition the Type 3 trapezoid $T_3$ into rectangles $C_0, \cdots, C_t$, $D_0, \cdots, D_t$ and $F_1$, triangles $E_0, \cdots, E_t$ with height $h_1=\lfloor \frac{x^{-1/6}}{\tan \theta} \rfloor$ and $F_2$ with height $h_2$ satisfying $0 \leq h_2<h_1$, as illustrated in Fig. \ref{f:6}. Here $t$ satisfies
$$
t < h/h_1= \frac{1}{2}x^{1/2}\Big/\left(\frac{x^{-1/6}}{\tan \theta}-r'\right)=\frac{x^{2/3}\tan \theta}{2(1-r'x^{1/6}\tan \theta)}
\leq \frac{1}{2}x^{2/3}\tan \theta,
$$
where $r'$ is the decimal part of $\frac{x^{-1/6}}{\tan \theta}$. The width of $C_k$, denoted by $c_k$, is set to be $\lfloor x^{1/3}+\sqrt{2}~x^{1/6}\rfloor-\lfloor x^{1/3}+(\sqrt{2}-k)x^{1/6}\rfloor$, and therefore $d_k$, the width of $D_k$, equals to $\lfloor x^{1/3}+(\sqrt{2}-k)x^{1/6} \rfloor+kh_1\tan\theta$, $k=0, \cdots, t$. Note that when $h_1>h$, then the number of $D_k$ is $0$, but the result still holds.

\begin{figure}[!htph]
\centering
\includegraphics[width=0.6\textwidth]{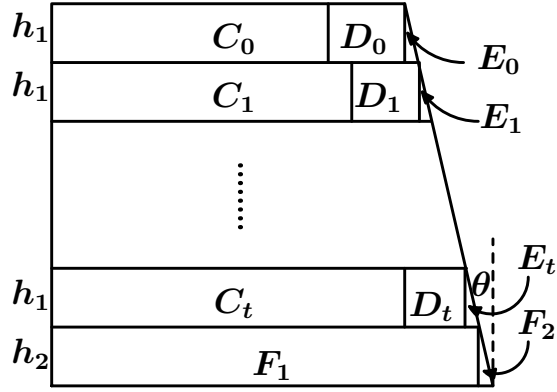}
\caption{Packing Type 3 trapezoid.}\label{f:6}
\end{figure}
1) Obviously, each $C_k$ can be packed perfectly with unit squares, thus
$$
W\left(\bigcup_{k=0}^{t} C_k\right)=\sum_{k=0}^{t} W(C_k)=0.
$$

2) It is easy to see that each $E_k$ can not be packed with unit squares. Thus
$$
W\left(\bigcup_{k=0}^{t} E_k\right)=\sum_{k=0}^{t} \frac{1}{2}h_1^2\tan \theta \leq  \frac{1}{4}x^{1/3}.
$$

3) We will estimate $W(\bigcup_{k=0}^{t} D_k)$ as follows. Since $d_0$ is an integer, $W(D_0)=0.$ For $k=1,\cdots,t$,  $0<kh_1\tan\theta<\frac{1}{2}x^{1/2}\tan \theta <1$ implies that $\lceil d_k \rceil=\lfloor x^{1/3}+(\sqrt{2}-k)x^{1/6} \rfloor+1$. Let $r_k$ be the decimal part of $x^{1/3}+(\sqrt{2}-k)x^{1/6}$. Then
\begin{equation}\label{e:3}
  \left\{ \begin {array}{lll}
d_k & = & x^{1/3}+(\sqrt{2}-k)x^{1/6}-r_k+kx^{-1/6}-kr'\tan\theta,\\
\lceil d_k \rceil & = & x^{1/3}+(\sqrt{2}-k)x^{1/6}-r_k+1.\\
\end {array}
\right.
\end{equation}
Next, we will pack $D_k$ with rectangles of size $1 \times \lceil d_k \rceil$ and estimate $\alpha_k$ more accurately than before. By \eqref{e:1}, we obtain
\begin{equation}\label{e:4}
\lceil d_k \rceil \cos \alpha_k+\sin \alpha_k=d_k.
\end{equation}
Substitute \eqref{e:4} into \eqref{e:3}, we have
\begin{equation}\label{e:5}
(x^{1/3}+(\sqrt{2}-k)x^{1/6}-r_k)(1-\cos \alpha_k)=\cos \alpha_k+\sin \alpha_k-kx^{-1/6}+kr'\tan\theta.
\end{equation}
Substitute Taylor's formulae for $\cos \alpha_k,\sin \alpha_k$,
$$
\left\{ \begin {array}{lll}
\cos \alpha_k &=& 1-\frac{1}{2}\alpha_k^2+\frac{1}{24}\alpha_k^4+o(\alpha_k^5),\\
\sin \alpha_k &=& \alpha_k-\frac{1}{6}\alpha_k^3+o(\alpha_k^4),\\
\end {array}
\right.
$$
into \eqref{e:5} and set $\alpha_k=l_{k1}x^{-1/6}+l_{k2}x^{-1/3}+l_{k3}x^{-1/2}+o(x^{-1/2})$. Since $0 < kr'\tan\theta <x^{-1/3}$, we set $kr'\tan\theta=\gamma_kx^{-1/3}+o(x^{-1/3})$, it follows that $0\leq \gamma_k<1$. Comparing the coefficients of terms $x^0$ and $
x^{-1/6}$, on both sides of \eqref{e:5}, we have
$$
\alpha_k=\sqrt{2}x^{-1/6}+0 \cdot x^{-1/3}+l_{k3}x^{-1/2}+o(x^{-1/2}).
$$
Since $0 \leq k < \frac{1}{2}x^{2/3}\tan \theta < \frac{\sqrt{2}}{2}x^{1/6}$, we set $ k=\beta_k x^{1/6}+o(x^{1/6})$, it follows that $0\leq \beta_k<\frac{\sqrt{2}}{2}$. Comparing the coefficients of terms $
x^{-1/3}$, on both sides of \eqref{e:5}, we have
$$
\alpha_k=\sqrt{2}x^{-1/6}+0 \cdot x^{-1/3}+\dfrac{r_k+\gamma_k-\frac{1}{6}\beta_k-\frac{5}{6}}{\sqrt{2}(1-\beta_k)}x^{-1/2}+o(x^{-1/2}).
$$
Hence $|\alpha_k-\alpha_{k-1}|\leq 3(1+\sqrt{2})x^{-1/2},k=2,\cdots,t$.

We pack $D_k$ as follows. First, we leave a Type 2 trapezoid $D_{11}$ at the top of $D_1$. Second, for $k=2,\cdots,t,$ we pack $D_{k-1}$ with rectangles of size $1 \times \lceil d_{k-1} \rceil$ when $b_k\geq \frac{1}{\cos \alpha_{k-1}}$. If not, we pack $D_k$ with rectangles of size $1 \times \lceil d_k \rceil$ (see Fig. \ref{f:7}). When $\alpha_{k-1} \geq \alpha_k$, the wasted region between $D_{k-1}$ and $D_k$ consists of a triangle $X_{k1}$ and trapezoids $X_{k2},X_{k3}$. The case of $\alpha_{k-1} < \alpha_k$ can be treated in similar fashion. Last, we leave Type 2 trapezoid $D_{t1}$ at the bottom of $D_t$.

\begin{figure}[!htph]
\centering
\includegraphics[width=0.5\textwidth]{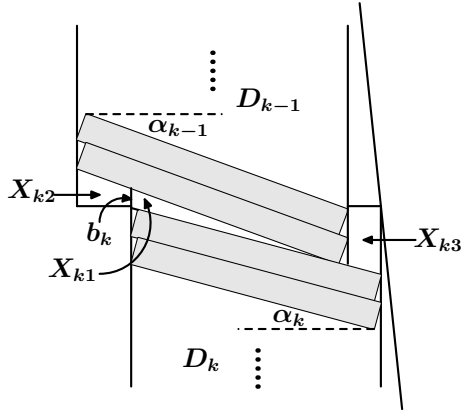}
\caption{The wasted region between $D_{k-1}$ and $D_k$.}\label{f:7}
\end{figure}

The total wasted area of both ends of rectangles of size $1 \times \lceil d_k \rceil,k=1,\cdots,t,$ is less than $\sum_{k=1}^{t} h_1\cdot2\cdot\frac{1}{2}\cdot1^2\tan\alpha_k <\frac{\sqrt{2}}{2}x^{1/3}$. By (ii), $W(D_{11})+W(D_{t1})\leq O(d_1^{5/6})+O(d_{t}^{5/6})=O(x^{5/18})$. The wasted area between
$D_{k-1}$ and $D_k$ is $S(X_{k1})+S(X_{k2})+S(X_{k3}) <
\frac{1}{2}(x^{1/3})^2 \cdot 3(1+\sqrt{2})x^{-1/2}+\frac{1}{2}(1+1+\sqrt{2})x^{1/6}+O(x^{1/6})O(x^{-1/6})\leq (\frac{5}{2}+2\sqrt{2})x^{1/6}$, which implies that the total wasted area of these joints is bounded by  $(\frac{5}{2}+2\sqrt{2})x^{1/6}\cdot t<(\frac{5}{4}\sqrt{2}+2)x^{1/3}$. Thus,
$$
W\left(\bigcup_{k=0}^{t} D_k\right)< 0+\frac{\sqrt{2}}{2}x^{1/3}+O(x^{5/18})+\left(\frac{5}{4}\sqrt{2}+2\right)x^{1/3}\leq \left(\frac{7}{4}\sqrt{2}+2\right)x^{1/3}.
$$

4) At last, we will estimate $W(F_1)$ and $W(F_2)$. The height of the rectangle $F_1$ satisfies $0 \leq h_2 <\min(h,h_1)$, and the width of it, denoted by $f_1$, satisfies $f_1 \sim x^{1/3}$. When $0 \leq h_2 \leq x^{1/3}$, we pack $\lfloor h_2 \rfloor \times \lfloor f_1 \rfloor$ unit squares into $F_1$, then $W(F_1)< h_2+f_1<2x^{1/3}$. When $x^{1/3}< h_2 \leq h$, we pack $F_1$ with rectangles of size $1 \times \lceil f_1 \rceil$, as shown in Fig. \ref{f:8}, where $F_{11},F_{12}$ are Type 2 trapezoids. Since $W(F_{11})+W(F_{12})=O(x^{5/18})$, the total wasted area of both ends of the rectangles of size $1 \times \lceil f_1 \rceil$ is less than $h\cdot \sqrt{2}x^{-1/6}\sim \frac{\sqrt{2}}{2}x^{1/3}$, so $W(F_1)<O(x^{5/18})+\frac{\sqrt{2}}{2}x^{1/3}<x^{1/3}$.
To sum up, $W(F_1)<2x^{1/3}$. We estimate $W(F_2)$ in two cases, too. When $0 <\theta < x^{-2/3}$,  $W(F_2)<S(F_2)<\frac{1}{2}h^2\tan\theta<\frac{1}{8}x^{1/3}$. When $x^{-2/3}\leq \theta < \sqrt{2} x^{-1/2}$, $W(F_2)<S(F_2)<\frac{1}{2}h_1^2\tan\theta<\frac{1}{2}x^{1/3}$. Therefore, $W(F_2)<\frac{1}{2}x^{1/3}$ which implies $W(F)\leq W(F_1)+W(F_2)<\frac{5}{2}x^{1/3}$.

\begin{figure}[!htph]
\centering
\includegraphics[width=0.7\textwidth]{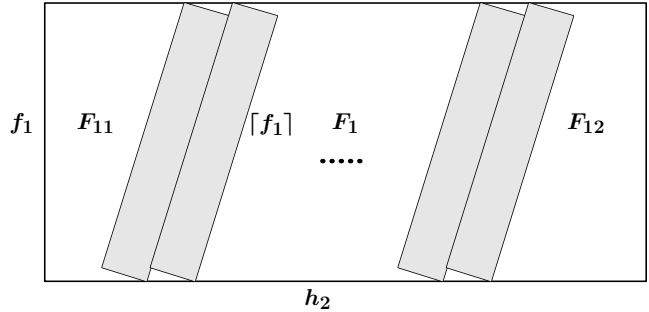}
\caption{Packing $F_1$ in the case of $x^{1/3}< h_2 \leq h $.}\label{f:8}
\end{figure}

Now, it follows from 1), 2), 3), 4) that the total wasted area
$$
W(T_3)\leq 0+\frac{1}{4}x^{1/3}+(\frac{7}{4}\sqrt{2}+2)x^{1/3}+\frac{5}{2}x^{1/3}
=(\frac{19}{4}+\frac{7}{4}\sqrt{2})x^{1/3}
$$
which completes the induction step. For $x \leq 100$, $W(T_1) \leq (1+c)x$. Because $c \leq 7,  (1+c)x^{3/8}<48<15\sqrt{2}+38<(15+c)\sqrt{2}+38, W(T_1) \leq (1+c)x<((15+c)\sqrt{2}+38)x^{5/8}$,
the proof of the initial step
 of the induction is completed.
\end{pf}

\section{Covering Waste Problem}\label{Covering}

Similarly, we can obtain the result of Covering Waste Problem. Note that in type 3 shape Trapezoid $T_3$, a top edge of length $a$ is modified, $a=\lfloor x^{1/3}-\sqrt{2}~x^{1/6}\rfloor$.

\begin{thm}\label{t:2}
Keep the notations above. Then
\begin{enumerate}[(i)]
  \item $W'(T_1)\leq ((15+c)\sqrt{2}+38)x^{5/8}$.
  \item $W'(T_2)\leq (\frac{19}{2}+\frac{7}{2}\sqrt{2})x^{5/6}$.
  \item $W'(T_3)\leq (\frac{19}{4}+\frac{7}{4}\sqrt{2})x^{1/3}$.
\end{enumerate}
Specially, when $T_1$ is a square of side length $x$, then $W'(x) \leq (16\sqrt{2}+38)x^{5/8}$.
\end{thm}

\begin{pf}  ~

(i) This can be proved in a similar argument to the one of (i) of Theorem \ref{t:1}.

(ii) This can be proved in a similar argument to the one of (ii) of Theorem \ref{t:1}.

(iii) We consider a coverage of Type 3 trapezoid $T_3$ with rectangles $C_k, D_k,k=1,\cdots,t,$ with height $h_1=\lfloor \frac{x^{-1/6}}{\tan \theta'} \rfloor$ and a rectangle $F_1$ with height $h_2$ satisfying $0 \leq h_2<h_1$. The width of $C_k$, denoted by $c_k$, is set to be $\lfloor x^{1/3}-\sqrt{2}~x^{1/6}\rfloor-\lfloor x^{1/3}-(\sqrt{2}+k)x^{1/6}\rfloor,$ and therefore
the width of $D_k$, denoted by $d_k$, is equal to $\lfloor x^{1/3}-(\sqrt{2}+k)x^{1/6} \rfloor+kh_1\tan\theta',k=1,\cdots,t$. It is easy to verify that the width of $F_1$, denoted by $f_1$, equals to $a+h\tan\theta'$ and $0 \leq t < \frac{1}{2}x^{2/3}\tan \theta'$. Set $E_k=D_k \setminus T_3,k=1,\cdots,t, F_2=F_1\setminus T_3$ (see Fig. \ref{f:9}), then
\begin{eqnarray*}
T_3 & = &\left(\bigcup_{k=1}^t C_k \right)\bigcup \left(\bigcup_{k=1}^t D_k \right )\bigcup F_1\setminus\left(\left(\bigcup_{k=1}^t E_k \right)\bigcup F_2\right),\\
W'(T_3) & \leq &\sum_{k=1}^t W'(C_k)+\sum_{k=1}^t S(E_k)+ W'\left(\bigcup_{k=1}^t D_k\right)+W'(F_1)+W'(F_2).
\end{eqnarray*}

\begin{figure}[!htph]
\centering
\includegraphics[width=.6\textwidth]{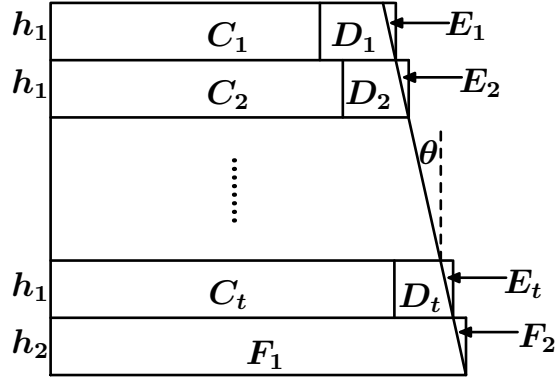}
\caption{Covering Type 3 trapezoid.}\label{f:9}
\end{figure}

1) Obviously, $\sum_{k=1}^t W'(C_k)=0$.

2) $\sum_{k=1}^t S(E_k)=\displaystyle\sum_{k=1}^{t} \frac{1}{2}h_1^2\tan \theta' \leq  \frac{1}{4}x^{1/3}$.

3) We estimate $W'(\bigcup_{k=1}^{t} D_k)$ as follows. For $k=1,\cdots,t$,
we want to cover $D_k$ with rectangles of size $1 \times \lceil d_k \rceil$ and estimate $\alpha_k$ more accurately. Similar to (iii) of Theorem \ref{t:1}, we can obtain
$$
|\alpha_k-\alpha_{k-1}|\leq 3(1+\sqrt{2})x^{-1/2}, k=2,\cdots,t.
$$

\begin{figure}[!htph]
\centering
\includegraphics[width=0.6\textwidth]{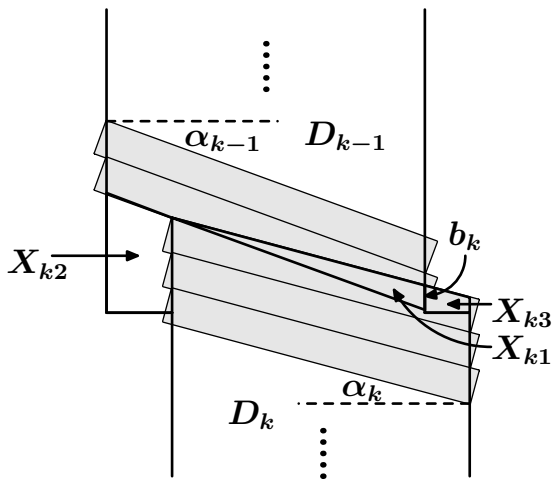}
\caption{The wasted area between $D_{k-1}$ and $D_k$ for covering.}\label{f:10}
\end{figure}

We cover $D_k$ as follows. First, we leave Type 2 trapezoid $D_{t1}$ at the bottom of $D_t$. Second, for $k=t,\cdots,2,$ we cover $D_k$ with rectangles of size $1 \times \lceil d_k \rceil$. When rectangles of size $1 \times \lceil d_k \rceil$ cover the right lower point of $D_{k-1}$, we cover $D_{k-1}$ with rectangles of size $1 \times \lceil d_{k-1} \rceil$, as shown in Fig. \ref{f:10}. When $\alpha_{k-1} \geq \alpha_k$, there are a triangle $X_{k1}$ and trapezoids $X_{k2}$, $X_{k3}$ between $D_{k-1}$ and $D_k$ needed to be solved further in the following. As shown in the figure, $b_k$ is the bottom edge of $X_{k3}$. The case of $\alpha_{k-1} < \alpha_k$ can be treated similarly. At last, we leave Type 2 trapezoid $D_{11}$ at the top of $D_1$.

The total wasted area of both ends of the rectangles of size $1 \times \lceil d_k \rceil,k=1,\cdots,t,$ is less than $\sum_{k=1}^{t} h_1\cdot2\cdot\frac{1}{2}\cdot1^2\tan\alpha_k <\frac{\sqrt{2}}{2}x^{1/3}$.
By (ii) of Theorem 2, $W'(D_{11})+W'(D_{t1})\leq O(d_1^{5/6})+O(d_{t}^{5/6})=O(x^{5/18})$. It is easy to see that $c_{k-1}-c_k$, the height of $X_{k2}$, is an exact integer. Let the bottom edge of $X_{k2}$ be $b'_{k2}$. We cover $X_{k2}$ with rectangles of size $(c_{k-1}-c_k) \times \lceil b'_{k2} \rceil$. The wasted area between
$D_{k-1}$ and $D_k$ is $S(X_{k1})+W'(X_{k2})+S(X_{k3}) <
\frac{1}{2}(x^{1/3})^2 \cdot 3(1+\sqrt{2})x^{-1/2}+\frac{1}{2}(1+1+\sqrt{2})x^{1/6}+O(x^{1/6})O(x^{-1/6})\leq (\frac{5}{2}+2\sqrt{2})x^{1/6}$, which implies that the total wasted area of these joints is bounded by $(\frac{5}{2}+2\sqrt{2})x^{1/6}\cdot t<(\frac{5}{4}\sqrt{2}+2)x^{1/3}$. Thus,
$$
W'\left(\bigcup_{k=0}^{t} D_k\right)< 0+\frac{\sqrt{2}}{2}x^{1/3}+O(x^{5/18})+\left(\frac{5}{4}\sqrt{2}+2\right)x^{1/3}\leq \left(\frac{7}{4}\sqrt{2}+2\right)x^{1/3}.
$$

4)At last, $W'(F_1)+W'(F_2)<\frac{5}{2}x^{1/3}$. The proof is similar to 4) of (iii) of Theorem \ref{t:1}.

By 1), 2), 3), 4), we obtain the total wasted area
$$
W'(T_3)\leq 0+\frac{1}{4}x^{1/3}+\left(\frac{7}{4}\sqrt{2}+2\right)x^{1/3}+\frac{5}{2}x^{1/3}
=\left(\frac{19}{4}+\frac{7}{4}\sqrt{2}\right)x^{1/3}.
$$
The proof of the induction step is omitted.
\end{pf}

\bibliography{mybibfile}

\end{document}